\documentclass[sn-mathphys-num]{sn-jnl}

\usepackage{mathtools}
\usepackage{graphicx}%
\usepackage{multirow}%
\usepackage{amsmath,amssymb,amsfonts}%
\usepackage{amsthm}%
\usepackage{mathrsfs}%
\usepackage[title]{appendix}%
\usepackage{xcolor}%
\usepackage{textcomp}%
\usepackage{manyfoot}%
\usepackage{booktabs}%
\usepackage{algorithm}%
\usepackage{algorithmicx}%
\usepackage{algpseudocode}%
\usepackage{listings}%
\usepackage{xspace}
\usepackage{rotating}
\usepackage{scalefnt}
\usepackage{upgreek}
\usepackage{MnSymbol}
\usepackage{enumitem}
\usepackage[utf8]{inputenc}
\definecolor{keywordcolor}{rgb}{0.7, 0.1, 0.1}   
\definecolor{tacticcolor}{rgb}{0.0, 0.1, 0.3}    
\definecolor{commentcolor}{rgb}{0.4, 0.4, 0.4}   
\definecolor{stringcolor}{rgb}{0.5, 0.3, 0.2}    
\definecolor{symbolcolor}{rgb}{0.1, 0.2, 0.7}    
\definecolor{sortcolor}{rgb}{0.1, 0.5, 0.1}      
\definecolor{attributecolor}{rgb}{0.7, 0.1, 0.1} 
\definecolor{errorcolor}{rgb}{1, 0, 0}           
\newcommand{\lean}[1]{\lstinline[language=lean]{#1}}
\newcommand{\mathlib}{\textsf{mathlib}\xspace}

\newtheorem{theorem}{Theorem}
\newtheorem{lemma}{Lemma}
\newtheorem{assumption}{Assumption}
\theoremstyle{thmstylethree}%
\newtheorem{definition}{Definition}%

\begin{document}

\title[Article Title]{Formalization~of~Complexity~Analysis~of~the~First-order~Algorithms~for~Convex~Optimization}


\author[1]{\fnm{Chenyi} \sur{Li}}\email{lichenyi@stu.pku.edu.cn}

\author[1]{\fnm{Ziyu} \sur{Wang}}\email{wangziyu-edu@stu.pku.edu.cn}

\author[2]{\fnm{Wanyi} \sur{He}}\email{2100017455@stu.pku.edu.cn}

\author[1]{\fnm{Yuxuan} \sur{Wu}}\email{2100010630@stu.pku.edu.cn}

\author[1]{\fnm{Shengyang} \sur{Xu}}\email{xushengyang0429@gmail.com}

\author*[3]{\fnm{Zaiwen} \sur{Wen}}\email{wenzw@pku.edu.cn}

\affil[1]{\orgdiv{School of Mathematical Sciences}, \orgname{Peking University}, \orgaddress{\country{China}}}
\affil[2]{\orgdiv{Yuanpei College}, \orgname{Peking University}, \orgaddress{\country{China}}}

\affil[3]{\orgdiv{Beijing International Center for Mathematical Research}, \orgname{Peking University}, \orgaddress{\country{China}}}


\abstract{The convergence rate of various first-order optimization algorithms is a pivotal concern within the 
numerical optimization community, as it directly reflects the efficiency of these algorithms across 
different optimization problems.  Our goal is to make a significant step forward in the formal mathematical representation 
of optimization techniques using the Lean4 theorem prover. We first formalize the gradient for smooth functions and the subgradient 
for convex functions on a Hilbert space, laying the groundwork for the accurate formalization 
of algorithmic structures. Then, we extend our contribution by proving several properties 
of differentiable convex functions that have not yet been formalized in Mathlib. Finally,  a comprehensive formalization of these algorithms is presented. These developments 
are not only noteworthy on their own but also serve as essential precursors to the formalization of 
a broader spectrum of numerical algorithms and their applications in machine learning as well as many other areas.\footnote{Our implementation of formalization of complexity analysis of the first-order algorithms for convex optimization can be found in \url{https://github.com/optsuite/optlib}}}




\maketitle

\section{Introduction}
Within the expansive domain of optimization and operational research, the analysis and application of first-order optimization algorithms are fundamental, crucial for addressing diverse challenges in fields such as machine learning \cite{bottou2018optimization}, data science, and engineering. 
Nevertheless, the theoretical foundations that ensure their efficacy, especially from the perspective of convergence analysis, are complex and require rigorous formalization. This paper is dedicated to navigating the complexities of formalizing the analysis of first-order optimization algorithms.
These algorithms are not merely tools for immediate problem-solving but also form the groundwork for developing more sophisticated optimization techniques.

To the best of our knowledge, few works relate to the formalization of convex optimization and numerical algorithms. However, the formalization of analysis has been extensively pursued by many researchers \cite{boldo2016formalization}
using various formalization languages, including Coq, Isabelle \cite{Nipkow2002APA} and Lean \cite{de2015lean}. For instance, Kudryashov formalized the divergence theorem and the Cauchy integral formula in Lean \cite{kudryashov2022formalizing}. Gouëzel extensively studied the formalization of the change of variables formula for integrals \cite{gouezel2022formalization}. The application of formal methods in machine learning was explored by Tassarotti \cite{tassarotti2021formal}. In the area of convex analysis and optimization, Grechuk presented a formalization of lower semicontinuous functions in Isabelle, including some related properties \cite{Lower_Semicontinuous-AFP}. Allamigeon provided a formalization of convex polyhedra based on the simplex method in Coq \cite{allamigeon2019formalization}. Verified reductions for optimization problems have also been explored \cite{verifiedoptimization}.

In this paper, building on Lean4 language and the corresponding \mathlib4 library \cite{mathlibcommunity}, we formalize the complexity analysis of first-order algorithms for convex and strongly convex functions, including the gradient descent method, the subgradient method, the proximal gradient method, and the Nesterov acceleration method \cite{nesterov1983method}. The theoretical properties of these numerical algorithms are discussed in various sources, including \cite{optimization_method, nesterov2018lectures, nocedal1999numerical, beck2017first, sun2006optimization, boyd2004convex}
. The main contributions of this paper are listed as follows.

1) To address derivative calculations in optimization, we propose formalizing the definition of the gradient. 
In \mathlib, the differentiability of a function is formalized using the \lean{fderiv} construct, which represents the function's derivative as a continuous linear map at differentiable points. However, the type-checking mechanisms inherent in Lean pose challenges for direct computation. This limitation underscores the need for a more computationally friendly representation of the gradient within \mathlib. By utilizing the Riesz Representation Theorem in a Hilbert space, we can transform the continuous linear map into a vector in the Hilbert space, thereby simplifying calculations with elements in this space.

2) We explore the formalization of the properties of convex functions and subgradients. The formalization of complexity analysis for first-order optimization algorithms fundamentally draws on the properties of convex functions. Currently, \mathlib's treatment of convex functions primarily encompasses their zero-order characteristics. This focus results in a notable absence of properties that leverage the function's gradient. Thus, we formalize properties such as the first-order conditions in this paper. Additionally, to address challenges associated with non-smooth optimization, we have extended the library by introducing the definitions of the subgradient and the proximal operator, alongside proofs of their relevant properties.

\begin{sidewaysfigure}
  \centering
  \includegraphics[width=\textwidth,height=0.6\textwidth]{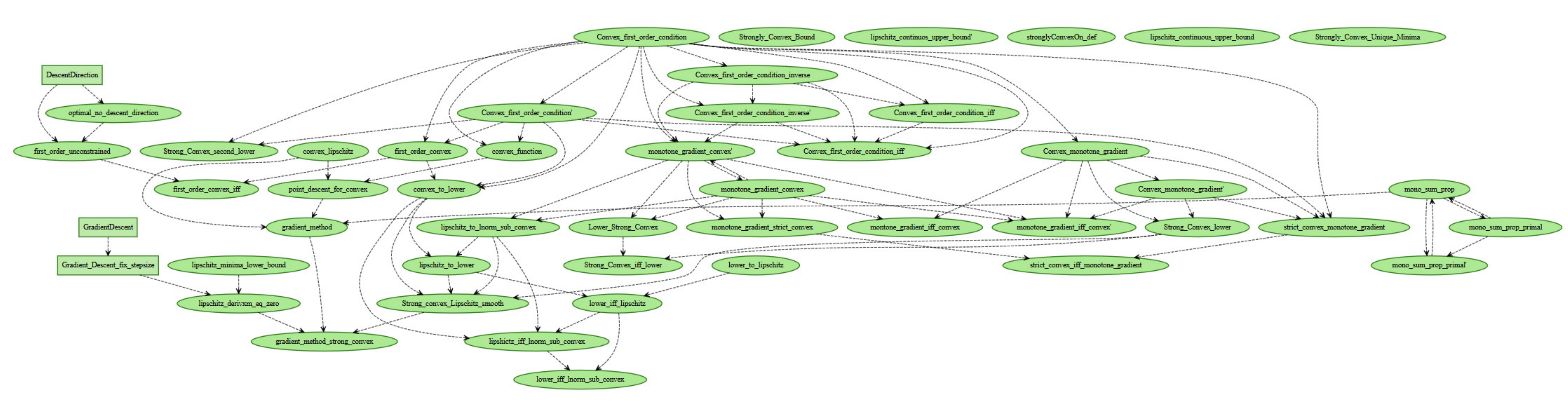}
  \caption{Part of the blueprint of the whole project (focus on the properties of convex functions and the convergence rate of the gradient descent method)}
  \label{fig: blueprint}
\end{sidewaysfigure}

3) Whereas the majority of current formalization efforts concentrate on theoretical mathematics, our work seeks to extend formalization into the realm of applied mathematics by formalizing numerical algorithms. This approach opens up broader possibilities for formalization across a wider range of fields. 
To broaden the applicability of our algorithm definitions to concrete optimization problems, we employ the \lean{class} structure to formalize the definitions of first-order algorithms, which facilitates a more generic representation. For implementing specific algorithm examples, the \lean{instance} structure \cite{baanen2022use} allows for the straightforward application of these algorithms, enabling users to instantiate specific cases and subsequently prove the requisite properties associated with them. We also build a blueprint for the whole project\footnote{The whole blueprint can be found in \url{https://chenyili0818.github.io/optlib-blueprint/dep_graph_document.html}}, which gives a brief introduction and contains the correlation between the definitions, theorems and proofs. Part of the blueprint, which focuses on the properties of convex functions and the convergence rate of the gradient descent method, is illustrated in Figure \ref{fig: blueprint}.     

The rest of the paper is organized as follows: In section \ref{sec: Mathematical preliminaries}, we briefly review the basic mathematical definitions and a general introduction to four types of first-order optimization algorithms. In section \ref{sec: Lean preliminaries}, we introduce relevant definitions already formalized by pioneers in the \mathlib community. The formalization of the definition and basic properties of the gradient and subgradient is presented in section \ref{sec: Gradient and Subgradient in Lean}. In sections \ref{sec: Properties of Convex Functions} and \ref{sec: Properties of Lipschitz Smooth Functions in Lean}, we formalize the properties of convex functions and L-smooth functions in Lean, respectively. The proximal operator is formally introduced in section \ref{sec: Proximal Operator in Lean}. Finally, in section \ref{sec: Convergence of First Algorithms in Lean}, we build the class for different first-order algorithms and prove the convergence rate of these algorithms.

\section{Mathematical preliminaries}
\label{sec: Mathematical preliminaries}

\subsection{The subgradient and the proximal operator}
Differentiability of a function within an Euclidean space is often characterized using little-o notation. When dealing with functions defined on general normed spaces, the complexity increases. To address this issue, we utilize the concept of the Fréchet derivative.

\begin{definition}
Let $E$ and $F$ be normed vector spaces, with $U \subseteq E$ representing an open subset. 
A function $f: U \to F$ is called Fréchet differentiable at a point $x \in U$ if there exists 
a bounded linear operator $A: E \to F$ satisfying the condition:
  \begin{align}
    \lim_{\|h\| \to 0} \frac{\|f(x+h)- f(x) -Ah\|_{F}}{\|h\|_E} = 0.
    \label{def : fderiv}
  \end{align}
\end{definition}

The concept of a subgradient is introduced to address points where the function may not be differentiable, yet still possesses certain advantageous properties.

\begin{definition}
For a function $f$ mapping from a Hilbert space $E$ to $\mathbb{R}$ and with $x$ in the domain of $f$, 
a vector $v$ is called a subgradient of $f$ at $x$ if for all $y \in E$,
  \begin{align*}
    f(y) \geq f(x) + \langle v, y-x \rangle.
  \end{align*}
\end{definition}

Define the collection of all subgradients at a point 
$x$ as the subderiv at that point, denoted
as $\partial f(x)$. It is critical to note that for convex functions, the subdgradient is guaranteed to be well-defined and nonempty at every point within the domain. Notably, at points where the function is smooth, the subderiv reduces to a singleton set containing only the gradient of the function at that point, i.e. $\partial f(x) = \{\nabla f (x)\}$. Building upon this conceptual framework, we next introduce the proximal operator.

\begin{definition}
For a function $f$ mapping from a Hilbert space to $\mathbb{R}$, the proximal operator is defined as:
  \begin{align*}
    \operatorname{prox}_f(x) = \underset{u}{\arg\min} \left\{f(u) + \frac{1}{2}\|u-x\|^2\right\}.
  \end{align*}
\end{definition}

For a convex function $f$, the addition of the $\|u-x\|^2$ term transforms the optimization problem into a strongly convex one, simplifying the original problem. Due to the characteristics of convex and strongly convex functions, the proximal operator is well-defined across all points, providing a means to minimize the function $f$ within a vicinity of the current point $x$. 

\subsection{First order algorithms solving optimization problems}
In this subsection, we give a brief review of the general first order algorithms solving optimization problem.
We mainly focus on unconstrained optimization problems
  \begin{align}
    \min f(x),
    \label{eq : unconstrained}
  \end{align}
where $f(x)$ is a convex function. Moreover, 
composite optimization problems are also considered:
  \begin{align}
    \min f(x) + g(x),
    \label{eq : composite}
  \end{align}
where $f(x)$ and $g(x)$ are convex, and $f(x)$ is smooth, while $g(x)$ may not be differentiable. The proximal gradient and Nesterov acceleration methods are particularly prevalent for these composite optimization problems. The efficiency of these algorithms, often measured by their convergence rates, is a key focus within the field of numerical optimization, making a detailed analysis of these rates essential.

\begin{itemize}
  \item \textbf{Gradient Descent Method} \\
  This foundational algorithm targets smooth functions in problem (\ref{eq : unconstrained}) and is notable for its simplicity and effectiveness. 
  The update mechanism is defined as
  \begin{align}
    x_{k+1} = x_k - \alpha_k \nabla f(x_k),
    \label{alg: gradient descent}
  \end{align}
  where $\alpha_k$ represents the stepsize for the $k$-th iteration, and $\nabla f(x_k)$ denotes the gradient at the point $x_k$. 
  Its convergence is characterized by $O(1/k)$ for convex functions and $O(\rho^k)$ for strongly convex functions, where $\rho$ indicates the condition number of the target function.
  \item \textbf{Subgradient Descent Method} \\
  In cases where the target function in problem (\ref{eq : unconstrained}) is nonsmooth and a gradient may not exist at every point, 
  the subgradient is utilized instead. The update formula is as follows
  \begin{align}
    x_{k+1} = x_k - \alpha_k g_k, \quad g_k \in \partial f(x_k),
    \label{alg: subgradient method}
  \end{align}
  where $g_k$ is the subgradient at $x_k$. The convergence rate for convex functions follows a $O(1/\sqrt{k})$ pattern. More concrete results can be found in \cite{boyd2003subgradient}.
  \item \textbf{Proximal Gradient Method} \\
  The proximal gradient method is widely used in optimization problems with the form (\ref{eq : composite}). 
  The update date scheme of this algorithm is given as
  \begin{align}
    x_{k+1} = \operatorname{prox}_{\alpha_k g}(x_k - \alpha_k \nabla f(x_k)),
    \label{alg: proximal}
  \end{align}
  where $\operatorname{prox}_{\alpha_k g}(x)$ denotes the proximal operator of the function $g$ at the point $x$. This 
  method can be viewed as an implicit version of subgradient method. The convergence rate of this algorithm is $O(1/k)$ under the assumptions stated above. More concrete results are referred to \cite{parikh2014proximal}.
  \item \textbf{Nesterov Acceleration Method} \\
  As an enhancement of the proximal gradient method, the Nesterov acceleration approach improves the convergence speed. 
  Nesterov acceleration method utilizes two sequences of points, $x_k$ and $y_k$, to update the point.
  The algorithm updates as following
  \begin{align}
    \begin{cases}   
      y_k = x_k + \frac{\gamma_k (1 - \gamma_{k - 1})}{\gamma_{k - 1}} (x_k - x_{k - 1}), \\
      x_{k} = \operatorname{prox}_{\alpha_k g}(y_k - \alpha_k \nabla f(y_k)).
    \end{cases}
    \label{alg: Nesterov}
  \end{align}
   Assuming the hyperparameters satisfy $\frac{(1 - \gamma_n)t_n}{\gamma_n^2} \leq \frac{t_{n - 1}}{\gamma_{n - 1}^2}$, the algorithm achieves the convergence rate of $O(\frac{1}{k^2})$. This method is also known as FISTA \cite{beck2009fast} which is widely used in compressive sensing. There is also another version of Nesterov acceleration scheme known as the second version of Nesterov acceleration, which is given as
  \begin{align}
  \label{alg: second Nesterov}
    \begin{cases}
      z_k = (1 -\gamma_k) x_{k-1} + \gamma_k y_{k-1},  \\
      y_{k} = \operatorname{prox}_{\frac{t_k}{\gamma_k} g}(y_{k-1} -\frac{t_k}{\gamma_k} \nabla f(z_k)),  \\
      x_{k} = (1 -\gamma_k) x_{k-1} + \gamma_k y_k.
    \end{cases}
  \end{align}
  The same convergence rate holds, if the hyperparameters satisfy
  $\frac{(1 - \gamma_{n + 1}) t_{n + 1}}{\gamma_{n + 1}^2} \leq \frac{t_n}{\gamma_n^2}$.
\end{itemize}

\section{Lean preliminaries}
\label{sec: Lean preliminaries}

\subsection{The differentiable structure of a normed space}

The \mathlib library, a comprehensive mathematical library for the Lean theorem prover, offers a robust framework for formalizing various concepts in calculus and analysis. Central to its calculus library is the concept of the Fréchet derivative, or fderiv, which facilitates the rigorous definition of the derivative for smooth functions between normed spaces.

In Lean, the fderiv structure is pivotal in defining the derivative of a smooth function between normed spaces. It encapsulates the derivative as a continuous linear map, adhering to the rigorous mathematical foundation for which Lean is renowned. The fderiv structure is defined as follows:

\begin{lstlisting}
structure HasFDerivAtFilter [NontriviallyNormedField k] [NormedAddCommGroup E] 
    [NormedSpace k E] [NormedAddCommGroup F] [NormedSpace k F] 
    (f: E → F) (f': E →L[k] F) (x: E) (L: Filter E): Prop where 
  of_isLittleO :: isLittleO: 
    (fun x' => f x' - f x - f' (x' - x)) =o[L] fun x' => x' - x
\end{lstlisting}

The utilization of a continuous linear map to define the derivative in Lean's 
\mathlib library enhances both generality and mathematical precision. Spaces \lean{E} and \lean{F}
are not limited to Euclidean spaces but can be any normed spaces over a nontrivially normed field $\mathbb{K}$. This broad applicability supports a wide range of mathematical and analytical discussions within the Lean environment. However, this generality introduces certain challenges in the context of numerical optimization. The abstract nature of continuous linear maps may lead to complications when devising update schemes for optimization algorithms. Precise type checks, a cornerstone of Lean's system, necessitate a reevaluation of the fderiv type when applied to numerical methods. 

Moreover, the \mathlib introduce the definition of \lean{deriv} to deal with the special case that \lean{E} is merely a NontriviallyNormedField $\mathbb{K}$. In this way, the continuous linear map becomes a single element in the space \lean{F}. 

To address these challenges, we pivot towards the gradient in vector form within \lean{E}. This approach aligns more closely with the practical requirements of numerical optimization, allowing for a more straightforward computation of update schemes. The transition from the Fréchet derivative to the gradient, along with the implications for numerical optimization, will be explored in detail in section \ref{sec: Gradient}.

\subsection{The convexity of a function}
The concept of convexity plays a pivotal role in optimization, underpinning many algorithms and theoretical results.

In the \mathlib library, the definition of a convex function is articulated through below:
\begin{lstlisting}
def ConvexOn (K: Type u_1) [OrderedSemiring K] [AddCommMonoid E] [OrderedAddCommMonoid β] [SMul K E] [SMul K β] (s: Set E) 
    (f: E → β): Prop :=
  Convex K s ∧ ∀ ⦃x⦄, x ∈ s → ∀ ⦃y⦄, y ∈ s → ∀ ⦃a b: K⦄, 0 ≤ a → 0 ≤ b → a + b = 1 → f (a • x + b • y) ≤ a • f x + b • f y
\end{lstlisting}

It is worth noting that the conditions on the input and output spaces are mild, which may not even require normed spaces. However, in this paper, we primarily focus on convex functions from a Hilbert space to $\mathbb{R}$, which 
is a special case of this definition as \lean{ConvexOn ℝ s f}.
The formalization of convexity within \mathlib provides a solid foundation for discussing and proving various properties of convex functions, particularly those that are differentiable.

While \mathlib's current formalization encompasses the core concept of convexity and some differentiable properties concerning only single-variable convex functions, there is ongoing work to enrich the library with additional properties related to differentiable convex functions of multiple variables, or more generally, on normed or Hilbert spaces. These properties are crucial for analyzing the behavior of optimization algorithms, especially in proving their convergence. The discussion of these extensions and their implications for algorithmic analysis will be elaborated upon in section \ref{sec: Properties of Convex Functions}.

\section{Gradient and Subgradient in Lean}
\label{sec: Gradient and Subgradient in Lean}
\subsection{Gradient}
\label{sec: Gradient}
The earlier discussion highlights that while \lean{fderiv} broadly defines derivatives within normed spaces, our interest in numerical optimization primarily lies with Hilbert spaces, which offer more intricate structures compared to normed spaces. Specifically, for functions mapping from a Hilbert space to the fields $\mathbb{R}$ or $\mathbb{C}$—collectively referred to as $\mathbb{K}$—the formal type of their Frechet derivative (fderiv) is denoted as \lean{E →L[K] K}. In the process of formalizing the gradient descent algorithm, the objective is to compute the update step, which involves applying the formula $x - \alpha \nabla f(x)$. This computation requires additive and scalar multiplicative operations between the point $x$ and its derivative $\nabla f(x)$.

However, using type of a continuous linear map from a Hilbert space to $\mathbb{K}$ does not directly support these operations. Consequently, converting the continuous linear map into a vector in the Hilbert space becomes crucial. This is precisely where the definition of the gradient becomes relevant and useful, as it is inherently designed to facilitate such operations by converting the abstract derivative into a tangible vector in the Hilbert space, thereby enabling the additive and scalar multiplicative operations necessary for the gradient descent update formula.
\begin{definition}
Let $E$ be a Hilbert space with inner product $\langle \cdot \, , \cdot \rangle$, while $U \subseteq E$ representing an open subset. 
  A function $f: U \to \mathbb{R}$ owns a gradient at a point $x \in U$ if there exists 
  a vector $g \in E$ satisfying the condition:
  \begin{align*}
    \lim_{\|h\| \to 0} \frac{f(x+h)- f(x) - \langle g, h\rangle}{\|h\|_E} = 0.
  \end{align*}
\end{definition}
Leveraging the definition of the Fréchet derivative, and utilizing the Riesz Representation Theorem on a Hilbert space, it becomes evident that the continuous linear operator $A$, integral to the formulation of the Fréchet derivative (see \ref{def : fderiv}), can be represented as
\begin{align*}
  A h = \langle g, h \rangle \quad \forall h \in E.
\end{align*}
In Lean, we can define the gradient as follows:

\begin{lstlisting}
def HasGradientAtFilter [RCLike K] [NormedAddCommGroup F] [InnerProductSpace K F] [CompleteSpace F] (f : F → K) (f' x : F) (L : Filter F) :=  HasFDerivAtFilter f (toDual K F f') x L
\end{lstlisting}

The segment \lean{toDual K F f'} is to convert an element from the space \lean{F} into an element within the dual space \lean{F →L[K] K}. This conversion is facilitated by a canonical mapping that links a Hilbert space to its corresponding dual space. Based on this definition, it enables the extension to more nuanced definitions such as \lean{HasGradientWithinAt}
and \lean{HasGradientAt}, which are more frequently used in the formalization of optimization algorithms.

It is crucial to distinguish between the spaces within which the gradient and the Fréchet derivative are defined. Specifically, the gradient is defined within a complete inner product space, and this specification is necessary to leverage the Riesz Representation Theorem. In contrast, the Fréchet derivative is applicable to a broader category of normed spaces. It is evident that for a general normed space, we cannot define the gradient as mentioned above.

\subsection{Subgradient}
To the best of our current knowledge, there has not yet been a formalization of the subgradient definition within the \mathlib library. Serving as an extension of the gradient, the subgradient concept accommodates non-smooth functions. The precise definition of the subgradient for a convex function is articulated as follows:

\begin{lstlisting}
def HasSubgradientWithinAt (f: E → ℝ) (g: E) (s: Set E) 
  (x: E): Prop := ∀ y ∈ s, f y ≥ f x + ⟪g, y - x⟫

def SubderivWithinAt (f: E → ℝ) (s: Set E) (x: E): Set E :=
  {g: E| HasSubgradientWithinAt f g s x}
\end{lstlisting}

A core theorem related to the subgradient is the existence of the subgradient at the interior point of the domain. For simplicity, we only consider the case when the function is convex. 

\begin{lstlisting}
theorem SubderivWithinAt.Nonempty [NormedAddCommGroup E] [InnerProductSpace ℝ E] [CompleteSpace E]{f : E → ℝ} (hf: ConvexOn ℝ s f) (hc: ContinuousOn f (interior s)): ∀ x ∈ interior s, (SubderivWithinAt f s x).Nonempty 
\end{lstlisting}

In this theorem, we assume that the function is continuous within the interior of the domain \lean{s}. This is a technical assumption, as only mild conditions are imposed on the space \lean{E}. However, if the input space \lean{E} is finite-dimensional, it is established that the convex function is continuous within the interior of the domain, or equivalently, any possible discontinuity of the convex function occurs only at the boundary points. In the proof of the theorem, a crucial element is a lemma stating the supporting hyperplane theorem. Viewed as a geometric version of the Hahn-Banach theorem, we utilize the theorem \lean{geometric_hahn_banach_open} in \mathlib, which asserts that given disjoint convex sets s, t, where s is open, there exists a continuous linear functional which separates them.

Another important aspect is the equivalence of the subgradient and the gradient at points where the function is smooth. This highlights that the subgradient is a more general definition of a gradient for non-smooth convex functions.
\begin{lstlisting}
theorem SubderivWithinAt_eq_gradient {f'x : E} (hx: x ∈ interior s) (hf: ConvexOn ℝ s f) (h: HasGradientAt f (f'x) x) :
  SubderivWithinAt f s x = {f'x} 
\end{lstlisting}
Furthermore, the computation of the subgradient for two convex functions holds significant importance. In this context, we refer to the Moreau-Rockafellar theorem, which is instrumental for subsequent proofs involving the proximal operator. The underlying intuition behind this theorem is direct, but needs a novel construction in the proof.
\begin{theorem}
  Assume $f_1$ and $f_2$ are two convex functions define on $E$, then we have for any $x \in E$
  \begin{align*}
    \partial f_1 (x) + \partial f_2 (x) = \partial (f_1 + f_2) (x).
  \end{align*}
\end{theorem}
The theorem is formalized as: 
\begin{lstlisting}
theorem SubderivAt.add {f₁ f₂: E → ℝ} (h₁: ConvexOn ℝ univ f₁) (h₂: ConvexOn ℝ univ f₂)(hcon: ContinuousOn f₁ univ):
  ∀ (x: E), SubderivAt f₁ x + SubderivAt f₂ x = SubderivAt (f₁ + f₂) x
\end{lstlisting}
We focus on proving a more stringent variant of the original Moreau-Rockafellar theorem, imposing stricter conditions on the convex function's domain. To simplify our analysis and avoid the complexities associated with the interior points of the function's domain, we assume that the function is convex across the entire space. A more comprehensive formulation of the theorem would necessitate exploring the continuity of convex functions within the interior of the domain, an endeavor we reserve for future investigation. Additionally, for general nonconvex functions, we can also define Fréchet differentiability as outlined in \cite{bolte2014proximal}. 

\section{Properties of Convex Functions in Lean}
\label{sec: Properties of Convex Functions}

Throughout the discussion from this section to the concluding section, we uniformly assume, except in specific cases, that the input space \lean{E} constitutes a Hilbert space, and $f$ represents a function mapping from \lean{E} to \lean{E}. Consequently, the gradient of $f$ functions as \lean{E → ℝ}. In certain scenarios, we will consider the domain of the function as a subset within \lean{E}, designated as \lean{s}. These parameters are specified as follows:

\begin{lstlisting}
variable [NormedAddCommGroup E] [InnerProductSpace ℝ E] 
variable [CompleteSpace E] 
variable {f: E → ℝ} {f': E → E} {s: Set E} 
\end{lstlisting} 

\subsection{General Properties}
For convex functions, certain properties are crucial for establishing the convergence of algorithms. These properties are encapsulated in the following theorem:
\begin{theorem}
Let $f$ be a smooth function defined on a convex set $s$. The statements below are equivalent:
\begin{enumerate}
    \item $f$ is convex on $s$.
    \item For all $x, y \in s$, the function satisfies the first-order condition: $f(y) \geq f(x) + \nabla f(x)^\top(y-x)$.
    \item For all $x, y \in s$, the gradient of $f$ is monotonic: $(\nabla f(x) - \nabla f(y))^\top(x-y) \geq 0$.
\end{enumerate}
\end{theorem}

This collection of theorems has been formalized in the \lean{Convex_Function.lean} file.
\begin{lstlisting}
theorem Convex_first_order_condition_iff' (h₁: Convex ℝ s) (h: ∀ x ∈ s, HasGradientAt f (f' x) x) : ConvexOn ℝ s f 
  ↔ ∀ x ∈ s → ∀ y ∈ s → f x + inner (f' x) (y - x) ≤ f y 
\end{lstlisting}
\begin{lstlisting}
theorem monotone_gradient_iff_convex' (h₁: Convex ℝ s) (hf: ∀ x ∈ s, HasGradientAt f (f' x) x) : ConvexOn ℝ s f 
  ↔ ∀ x ∈ s, ∀ y ∈ s, inner (f' x - f' y) (x - y) ≥ (0: ℝ)
\end{lstlisting}

In these theorems, it is important to note that the use of the \lean{gradient} definition is not strictly necessary, as the term $\nabla f(x)^\top(y-x)$ is interpreted as the continuous linear map at $x$, evaluated at $y-x$, producing a real number. To provide a comprehensive formalization, we present statements for each theorem in both \lean{fderiv} and \lean{gradient} forms. For simplicity, we have shown only the version utilizing the gradient above. These theorems introduce a practical method for assessing the convexity of a function through gradient information. More automated ways of checking the convexity of a function can be explored in future work.

\subsection{Strongly Convex Functions}
While gradient descent exhibits a sublinear convergence rate for convex functions, it can achieve a linear convergence rate for strongly convex functions. The formalization of strongly convex functions represents a pivotal advancement in accurately formalizing the convergence rates of gradient descent across various function types. The definitions for uniform convexity and strong convexity are delineated as follows:
\begin{lstlisting}
def UniformConvexOn (s: Set E) (φ: ℝ → ℝ) (f: E → ℝ): Prop :=
  Convex ℝ s ∧ ∀ ⦃x⦄, x ∈ s → ∀ ⦃y⦄, y ∈ s → ∀ ⦃a b: ℝ⦄, 0 ≤ a → 0 ≤ b  → a + b = 1 → f (a • x + b • y) ≤ a • f x + b • f y - a * b * φ ‖ x - y ‖
\end{lstlisting}
\begin{lstlisting}
def StrongConvexOn (s: Set E) (m: ℝ): (E → ℝ) → Prop :=
  UniformConvexOn s fun r ↦ m / (2: ℝ) * r ^ 2
\end{lstlisting}
It is essential to clarify that the concept of uniform convexity can be applied within the framework of general normed spaces. However, strong convexity necessitates a definition within a Hilbert space, primarily due to the need to utilize the inner product to decompose the expression $\|x-y\|^2$. Following the establishment of this definition, it is imperative to elucidate the properties of strongly convex functions, leveraging derivative information. Consequently, we can formalize the following theorem concerning strongly convex functions

\begin{theorem}
Let $f$ be a function defined on a convex set $s$. The following statements are equivalent:
\begin{enumerate}
\item $f$ exhibits $m$-strong convexity on $s$.
\item The function $g(x) = f(x) - \frac{m}{2}\|x\|^2$ is convex on $s$.
\item For differentiable $f$, for all $x, y \in s$, it holds that $f(y) \geq f(x) + \nabla f(x)^\top(y-x) + \frac{m}{2} \|y - x\|^2$.
\item For differentiable $f$, for all $x, y \in s$, it holds that $(\nabla f(x) - \nabla f(y))^\top(x- y) \geq m \|x - y\|^2$.
\end{enumerate}
\end{theorem}
We only list the most important part of the formalization of the theorem here, while more detailed descriptions can be found in the Lean file \lean{Strong_Convex.lean}.
\begin{lstlisting}
theorem Strong_Convex_iff_lower (hf: ∀ x ∈ s, HasGradientAt f (f' x) x) (hs: Convex ℝ s) : StrongConvexOn s m f 
  ↔ ∀ x ∈ s, ∀ y ∈ s, inner (f' x - f' y) (x - y) ≥ m * ‖ x - y ‖ ^ 2 
\end{lstlisting} 

\section{Properties of Lipschitz Smooth Functions in Lean}
\label{sec: Properties of Lipschitz Smooth Functions in Lean}
Another significant class of function is the Lipschitz smooth function. The concept of Lipschitz smoothness serves to quantify a function's degree of smoothness. This property is formalized through the notion of Lipschitz continuity for a function over a specific set, which is defined in the \mathlib library as follows:
\begin{lstlisting}
def LipschitzOnWith [PseudoEMetricSpace α] [PseudoEMetricSpace β] (K: ℝ≥0) (f: α → β) (s: Set α) :=
  ∀ ⦃x⦄, x ∈ s → ∀ ⦃y⦄, y ∈ s → edist (f x) (f y) ≤ K * edist x y
\end{lstlisting}

A central theorem regarding Lipschitz smooth functions pertains to their upper bound. The lemma is articulated as follows:
\begin{lemma}
  Let $f$ be a $l$-Lipschitz smooth function defined on a  set $s$, then it holds
  \[f(y) \leq f(x) + \nabla f(x) ^ T(y-x) + \frac{l}{2} \lVert y - x \rVert^2, \quad \forall x, y \in s.\]
\end{lemma}
Within the formalized framework, we provide both the Frechet derivative and the gradient formulations of this theorem. For the sake of brevity, here we present only the fderiv formulation as:
\begin{lstlisting}
theorem lipschitz_continuous_upper_bound [NormedAddCommGroup E] [NormedSpace ℝ E] {f: E → ℝ} {f': E → E →L[ℝ] ℝ} {l: NNReal} (h₁: ∀ x₁: E, HasFDerivAt f (f' x₁) x₁) (h₂: LipschitzWith l f') :
  ∀ (x y: E), f y ≤ f x + (f' x) (y - x) + l / 2 * ‖ y - x ‖ ^ 2
\end{lstlisting}
In this proof, we use the auxiliary function $g(t) = f(x + t(y-x))$ as a function from $\mathbb{R}$ to 
$\mathbb{R}$. Using this function, we can transform the original problem to a one-variable problem,
and then utilize the mean-value theorem \lean{image_le_of_deriv_right_le_deriv_boundary} to get the result.

When it comes to convex Lipschitz smooth function, we can derive more properties of the function considering the
convexity of the function. We state the theorem as:
\begin{theorem}
  Let $f$ be a differentiable convex function defined on $\mathbb{R}^n$, then the following statement is equivalent
    \begin{enumerate}
      \item $f$ is $l$ - Lipschitz smooth on $\mathbb{R}^n$.
      \item $g(x) = \frac{l}{2}\|x\|^2 - f(x) $ is convex .
      \item $(\nabla f(x) - \nabla f(y)) ^ T(x- y) \geq \frac{1}{l} \lVert \nabla f(x) - \nabla f(y) \rVert^2, \quad \forall x, y \in \mathbb{R}^n$.
    \end{enumerate}
\end{theorem}
Note that sometimes the natural language statement would hide some of the assumptions which human 
would think as trivial, but in formalization, such assumptions need to be stated explicitly. We can state the 
formalization of the above theorem as :
\begin{lstlisting}
variable [ProperSpace E] {l: ℝ}
\end{lstlisting}
\begin{lstlisting}
theorem lower_iff_lipschitz (h₁: ∀ x, HasGradientAt f (f' x) x) (hfun: ConvexOn ℝ Set.univ f) (hl: l > 0): LipschitzWith l f' 
  ↔ ∀ x y, inner (f' x - f' y) (x - y) ≥ 1 / l * ‖ f' x - f' y ‖ ^ 2
\end{lstlisting}
\begin{lstlisting}
theorem lipshictz_iff_lnorm_sub_convex (h₁: ∀ x, HasGradientAt f (f' x) x) (hfun: ConvexOn ℝ Set.univ f) (hl: l > 0) :
  LipschitzWith l f' ↔ ConvexOn ℝ univ (fun x ↦ l / 2 * ‖ x ‖ ^ 2 - f x) 
\end{lstlisting} 
For functions that are both strongly convex and have a Lipschitz continuous gradient, we can propose an enhanced estimation, specifically formulated in the following theorem:
\begin{lemma}
  Let $f$ be a $l$-Lipschitz smooth and $m$-strongly convex function defined on $\mathbb{R}^n$, then the following 
  inequality holds,
  \[
    (\nabla f(x) - \nabla f(y)) ^ T(x- y) \geq\frac{ml}{m+l}\|x-y\|^2+ \frac{1}{m + l} \lVert \nabla f(x) - \nabla f(y) \rVert^2, \quad \forall x, y \in \mathbb{R}^n.
  \]
\end{lemma}

The formalized theorem is stated as:
\begin{lstlisting}
theorem Strong_convex_Lipschitz_smooth 
    (hsc: StrongConvexOn univ m f) (hf: ∀ x, HasGradientAt f (f' x) x) (h₂: LipschitzWith l f') (hl: l > (0: ℝ)):
  inner (f' x - f' y) (x - y) ≥ m * l / (m + l) * ‖ x - y ‖ ^ 2 + 1 / (m + l) * ‖ f' x - f' y ‖ ^ 2 
\end{lstlisting}

\section{Proximal Operator in Lean}
\label{sec: Proximal Operator in Lean}
In this section, we need to introduce an additional assumption on the space \lean{E}, specifically \lean{[CompleteSpace E]}. The rationale behind this will be explained later. To define the proximal operator in Lean, we must take a few steps to circumvent the direct use of $\operatorname{argmin}$ as commonly described in natural language. Since the operator $\operatorname{argmin}$ must be clarified as to whether the target function can reach the minima at a finite point, defining it directly is not straightforward. Instead, we can define the proximal property and then identify the set of points that satisfy this property. If we can demonstrate that this set is non-empty, we can then select one of these points as the proximal point.

Firstly, we can define the proximal property as:
\begin{lstlisting}
def prox_prop (f: E → ℝ) (x: E) (xm: E): Prop :=
  IsMinOn (fun u ↦ f u + ‖ u - x ‖ ^ 2 / 2) univ xm
\end{lstlisting}
We define the proximal set as all the points that satisfy the proximal property. This set is unique when the function $f$ possesses certain desirable properties, and it may be empty when $f$ is neither continuous nor convex.
\begin{lstlisting}
def prox_set (f: E → ℝ) (x: E): Set E := {u | prox_prop f x u}
\end{lstlisting}
For the proximal point, assuming that the proximal set is nonempty, we simply need to select one of its members. Here, we use the function \lean{Classical.choose} to select one element from this nonempty set.

\begin{lstlisting}
def prox_point (f: E → ℝ) (x: E) (h: ∃ y, prox_set f x y): E :=
  Classical.choose h
\end{lstlisting}

After defining the proximal operator, we need to prove the wellposedness of the proximal operator. Generally speaking, we have the following theorem.
\begin{theorem}
   The proximal set of each point is nonempty and compact, when $f$ satisfies one of the following conditions:
\begin{enumerate}
    \item $f$ is lower semicontinuous and has lower bound over the whole space.
    \item $f$ is a continuous convex function, in this case, the proximal set is unique.
  \end{enumerate}
\end{theorem}
In this theorem, we evaluate whether the function's minima can be achieved. For this reason, we need to consider the relationship between closed and compact sets, specifically, the requirement that any bounded closed set be compact. This condition is straightforward in Euclidean space but does not generally hold in infinite-dimensional Hilbert spaces. This equivalence implies that the Hilbert space is finite-dimensional. In \mathlib, we utilize the definition of a "proper space," which is characterized as a space in which every bounded closed set is compact. It is evident that some Hilbert spaces, such as $L^2([0,1])$, are not proper spaces, whereas Euclidean space is an example of a proper space.
When formalizing, it becomes necessary to relax the conditions required by the theorem, as we are not working within $\mathbb{R}^n$ but in a more abstract space $E$. This adjustment allows us to appreciate the properties $\mathbb{R}^n$ inherently possesses. We can then articulate the following formalized theorem
\begin{lstlisting}
theorem prox_set_compact_of_lowersemi (f: E → ℝ) (hc: LowerSemicontinuous f) (lbdf: BddBelow (f '' univ)) : ∀ x, Nonempty (prox_set f x) ∧ IsCompact (prox_set f x) 
\end{lstlisting} 
\begin{lstlisting}
theorem prox_set_compact_of_convex (f: E → ℝ) (hc: ContinuousOn f univ) (hconv: ConvexOn ℝ univ f) :
  ∀ x, Nonempty (prox_set f x) ∧ IsCompact (prox_set f x) 

theorem prox_unique_of_convex (f: E → ℝ) (x: E) (hfun: ConvexOn ℝ univ f)
  (h₁: prox_prop f x y₁) (h₂: prox_prop f x y₂): y₁ = y₂
\end{lstlisting} 

We can derive additional properties of the proximal operator, particularly its connection to the subgradient. This relationship is readily apparent from the optimality conditions of the unconstrained optimization problem.
\begin{theorem}
  If $f$ is a closed and convex function, then we have  
  \begin{align*}
    u = \operatorname{prox}_f(x) \quad \Leftrightarrow \quad x - u \in \partial f (u).
  \end{align*}
\end{theorem}
In Lean, we state the theorem as:
\begin{lstlisting}
theorem prox_iff_subderiv (f: E → ℝ) (hfun: ConvexOn ℝ univ f) : ∀ u: E, prox_prop f x u ↔ x - u ∈ SubderivAt f u := by
\end{lstlisting}

\section{Convergence of First Order Algorithms in Lean}
\label{sec: Convergence of First Algorithms in Lean}
In this section, we  give the formalization of first order algorithms in Lean using \lean{class} structure.
From the perspective that the \lean{class} structure in Lean is easy to generalize for different target functions.
For specialized problem, such as  the LASSO problem in compressive sensing with target function 
$f(x) = \|Ax -b\|^2$ and $g(x) = \|x\|_1$, we can use the \lean{instance} structure to formalize 
the algorithm for this specific problem. For each algorithm, under different assumptions on the stepsize,
we will get the convergence theorem. In this section, we assume that $E$ is a Hilbert space and $f$ is a 
function defined on $E$. $xm$ is a point in $E$ and denotes for the minima of the function,
and $x_0$ denotes the initial point we put into the algorithm. Generally speaking, an algorithm contains the following parts.
\begin{itemize}
  \item \textbf{Update scheme}: we will take track of the update points in the algorithm.
  \item \textbf{Information on the target function}: we need information for the target function, like the gradient, 
      and the Lipschitz continuous information on the gradient.
  \item \textbf{Step size constraint}: only suitable stepsize choice is admittable for the corresponding algorithm.
\end{itemize}
\subsection{Gradient Descent Method}
For general definition of gradient descent method in Lean, we use the \lean{class} type to define what 
a numerical optimization method is in Lean. In this class we have the function $f$, the gradient $f'$, 
and the initial point as the input, and contains the necessary information with the optimization problem.

\begin{lstlisting}
class GradientDescent (f : E → ℝ) (f' : E → E) (x0 : E) := 
  (x : ℕ → E) (a : ℕ → ℝ) (l : NNReal) 
  (diff : ∀ x₁, HasGradientAt f (f' x₁) x₁) 
  (smooth : LipschitzWith l f') (hl : l > 0)
  (update : ∀ k : ℕ, x (k + 1) = x k - a k • f' (x k)) 
  (step₁ : ∀ k, a k > 0) (initial : x 0 = x0)
\end{lstlisting}
\begin{lstlisting}
class GD_fixed_stepsize (f : E → ℝ) (f' : E → E) (x0 : E) :=
  (x : ℕ → E) (a : ℝ) (l : NNReal) (initial : x 0 = x0)
  (diff : ∀ x₁, HasGradientAt f (f' x₁) x₁) 
  (smooth : LipschitzWith l f') (hl : l > (0 : ℝ))
  (update : ∀ k : ℕ, x (k + 1) = x k - a • f' (x k)) 
  (step₁ : a > 0) (initial : x 0 = x0)
\end{lstlisting}
We can also define the gradient descent with fixed step size as a special instance of the general gradient
descent method. In this paper, we mainly focus on the fixed step size version of the gradient descent method,
but more general version can be added easily based on this work.

\begin{lstlisting}
instance {f: E → ℝ} {f': E → E} {x₀: E} 
    [p: GD_fixed_stepsize f f' x₀] : GradientDescent f f' x₀
\end{lstlisting}
It is straightforward to see that the gradient descent method with fixed stepsize is a special case of 
the gradient descent method, hence we can get the \lean{instance} structure as above.

The convergence rate of the fixed step size gradient method is give by the following theorem:
\begin{theorem}
  For unconstrained optimization problem (\ref{eq : unconstrained}), $f(x)$ is L-smooth. Let $f^* = f(x^*) = \inf_x f(x)$ be the optimal value. 
  \begin{enumerate}
    \item If $f(x)$ is convex, then for any step size $\alpha$ satisfying $0 < \alpha \leq \frac{1}{L}$, 
        the gradient descent algorithm (\ref{alg: gradient descent}) generates a sequence of points $\{x_k\}$ whose  function values satisfy the inequality
        \begin{align*}
          f(x_k) - f (x^*) \leq \frac{1}{2\alpha k} \|x_0 - x^*\|^2, \quad \forall k \in \mathbb{N}.
        \end{align*}
    \item If $f(x)$ is $m$-strongly convex, then for any step size $\alpha$ satisfying $0 < \alpha \leq \frac{2}{m + L}$, 
        the gradient descent algorithm (\ref{alg: gradient descent}) generates a sequence of points $\{x_k\}$ whose function values satisfy the inequality
        \begin{align*}
          \|x_k - x^*\|^2 \leq \left(1 - \alpha \frac{2mL}{m + L}\right)^k \|x_0 - x^*\|^2, \quad \forall k \in \mathbb{N}.
        \end{align*}
  \end{enumerate}
  
\end{theorem}
To prove the convergence rate of the fixed step size gradient descent method, we need to prepare for a bunch of 
theorems ahead, including the one iteration property of the method and the sum up property of the monotonic sequences.
Finally we can prove the convergence rate of the gradient descent method for Lipschitz smooth function. 
\begin{lstlisting}
theorem gradient_method_fix_stepsize 
    {alg: GD_fixed_stepsize f f' x₀} 
    (hfun: ConvexOn ℝ Set.univ f) (step₂: alg.a ≤ 1 / alg.l) :
  ∀ k: ℕ, f (alg.x (k + 1)) - f xm ≤ 1 / (2 * (k + 1) * alg.a) 
  * ‖ x₀ - xm ‖ ^ 2 
\end{lstlisting}
It is interesting to find out from the proof that there is no assumptions on \lean{xm} here. In general setting,
we use the case which \lean{xm} is the minima, but in the proof, we can see that the proof is valid for any point \lean{xm}.
So doing formalized proof can let us know the direct connection between the assumptions and the theorem.
\begin{lstlisting}
theorem gradient_method_strong_convex_fix_stepsize 
    {alg: GD_fixed_stepsize f f' x₀} 
    (hsc: StrongConvexOn univ m f) (hm: m > 0) 
    (min: IsMinOn f univ xm) (step₂: alg.a ≤ 2 / (m + alg.l)) : 
  ∀ k: ℕ , ‖ alg.x k - xm ‖ ^ 2 ≤ (1 - alg.a * (2 * m * alg.l / (m + alg.l))) ^ k * ‖ x₀ - xm ‖ ^ 2 
\end{lstlisting}
\subsection{Subgradient Descent Method}
In this subsection, we focus on the subgradient descent method. For subgradient descent method, our assumption is given as:
\begin{assumption}
Considering the unconstrained optimization problem (\ref{eq : unconstrained}), we assume
\label{assump: subgradient}
    \begin{enumerate}
  \item $f$ is convex on $s$.
  \item there exists at least one minima $x^*$ and $f(x^*) > - \infty$.
  \item $f$ is Lipschitz continuous, i.e. $|f(x) - f(y)| \leq G \|x - y\|$ for all $x, y \in \mathbb{R}^n$ with $G>0$. 
\end{enumerate}
\end{assumption}
Note that the assumption (c) is equivalent with assuming that the subgradient of the target function $f$ 
is bounded by $G$. The subgradient descent method is defined as follows:
\begin{lstlisting}
class subgradient_method (f: E → ℝ) (x₀: E) :=
  (x g: ℕ → E) (a: ℕ → ℝ) (ha: ∀ n, a n > 0)
  (G: NNReal) (lipschitz: LipschitzWith G f) (initial: x 0 = x₀)
  (update: ∀ k, (x (k + 1)) = x k - a k • (g k))
  (hg: ∀ n, g n ∈ SubderivAt f (x n))
\end{lstlisting}

Many different results can be derived with different kinds of step sizes. For simplicity, we only show theorem for the diminishing step size  in this paper, while more relevant results such as the convergence rate of fixed step size can be found in the code.
\begin{theorem}
    \label{thm: subgradient method}
  Suppose that Assumption \ref{assump: subgradient} is satisfied, and the step size sequence ${\alpha_k > 0}$ for all $k$,  $\alpha_k \to 0$  and $\sum_{k=0}^{\infty} \alpha_k = +\infty$, then the sequence $\{x^k\}$ generated by subgradient method (\ref{alg: subgradient method}) converges to the optimal solution $x^*$, and for all $k \geq 0$ with the rate 
  \begin{equation}
       \hat f^k - f^* \leq \frac{\| x^0 - x^* \|^2 + G^2 \sum_{i=0}^{k} \alpha_i^2}{2\sum_{i=0}^{k} \alpha_i},
      \end{equation}
    where $f^* = f(x^*)$, and $f^k$ is the minimum value of $f(x)$ up to the $k^{th}$ iteration of $f(x)$'s values, i.e., $\hat f^k = \min_{0 \leq i \leq k} f(x^i)$.
\end{theorem}


Then we have the formalized version of theorem (\ref{thm: subgradient method}) as: 
\begin{lstlisting}
theorem subgradient_method_converge :
  ∀ k, 2 * ((Finset.range (k + 1)).sum alg.a) *
  (sInf {x | ∃ i ∈ Finset.range (k + 1), f (alg.x i) = x} - f xm)
  ≤ ‖ x₀ - xm ‖ ^ 2 + alg.G ^ 2 * (Finset.range (k + 1)).sum (fun i => alg.a i ^ 2) 
\end{lstlisting}
Moreover we can have the convergence result as:
\begin{lstlisting}
theorem subgradient_method_diminishing_step_size
    (ha' : Tendsto alg.a atTop (nhds 0))
    (ha'' : Tendsto (fun (k : ℕ) => (Finset.range (k + 1)).sum alg.a) atTop atTop) :
  Tendsto (fun k => sInf {f (alg.x i) | i ∈ Finset.range (k + 1)}) atTop (nhds (f xm)) 
\end{lstlisting}
\subsection{Proximal Gradient Method}
\label{sec: proximal gradient method}
From this subsection to the end of the paper, for the usage of the proximal operator, we need to require the space \lean{E} satisfying \lean{[ProperSpace E]}. Considering the composite optimization problem (\ref{eq : composite}). Using the proximal property we defined and formalized in section \ref{sec: Proximal Operator in Lean}, 
we can give the formalization of the proximal gradient method (\ref{alg: proximal}) in Lean. In this method, we use the definition
\lean{prox_prop} rather \lean{prox_point} since for general function, the proximal set is not unique.
We have that any point in the proximal set satisfying the proximal property is admittable for proximal 
gradient method. Similar to the first order algorithms above, we also define a class for proximal 
gradient method.

\begin{lstlisting}
class proximal_gradient_method (f h : E → ℝ) (f': E → E) (x₀ : E) :=
  (xm : E) (t : ℝ) (x : ℕ → E) (L : NNReal)
  (fconv : ConvexOn ℝ univ f) (hconv : ConvexOn ℝ univ h) 
  (h₁ : ∀ x₁, HasGradientAt f (f' x₁) x₁) (h₂ : LipschitzWith L f') 
  (h₃ : ContinuousOn h univ) (minphi: IsMinOn (f + h) Set.univ xm) 
  (tpos: 0 < t) (step: t ≤ 1 / L) (ori: x 0 = x₀) (hL: L > (0: ℝ)) 
  (update: ∀ k, prox_prop (t • h) (x k - t • f' (x k)) (x (k + 1)))
\end{lstlisting}
First we need to give the basic assumptions for this problem.
\begin{assumption}
  \label{assumption}
  For composite optimization problem (\ref{eq : composite}), we have assumptions below:
  \begin{enumerate}
    \item $f$ is a differentiable convex function with $L$-Lipschitz continuous gradient.
    
    \item The function $h$ is continuous convex function (which means the proximal operator is well defined here);
    
    \item The minima of function $\phi$ is attainable at the finite point $x^*$, with the minimal value $\phi(x^*) > -\infty$.
  \end{enumerate}
\end{assumption}

We can get the convergence rate for proximal gradient as the theorem below:
\begin{theorem}
Suppose that Assumption \ref{assumption} is satisfied and the fixed step size 
$t_k = t  \in (0, \frac{1}{L}]$, then the sequence $\{x^k\} $ generated by (\ref{alg: proximal}) 
satisfies
\begin{align*}
  \psi(x^k) - \psi^* \leq \frac{1}{2kt} \| x^0 - x^* \|^2.
\end{align*}
\end{theorem}
The formalized convergence rate is given as:
\begin{lstlisting}
theorem proximal_gradient_method_converge 
    {alg: proximal_gradient_method f h f' x₀} :
  ∀ (k: ℕ+), (f (alg.x k) + h (alg.x k) - f alg.xm - h alg.xm) ≤ 1 / (2 * k * alg.t) * ‖ x₀ - alg.x ‖ ^ 2
\end{lstlisting}
\subsection{Nesterov Acceleration Method}
In this section, we mainly focus on the formalization of the Nesterov's acceleration method used on composite optimization. Since there are a few forms of the Nesterov's acceleration method,  we only choose two of them which 
are formalized in two relevant files \lean{Nesterov_Acceleration_first.lean} and \lean{Nesterov_Acceleration_second.lean}. Although having differences in the update scheme, they enjoy the same acceleration convergence rate. 

In this paper, we also exploit the \lean{instance} structure to define the algorithms. Firstly, we define the general method with abstract hyperparameter $\gamma$ and stepsize $t$, and then use the instance structure to connect the definition with the fixed stepsize ones. For the first form of the Nesterov's acceleration method, which is also known as FISTA method, we can formalize the fix stepsize version of the algorithm as:

\begin{lstlisting}
class Nesterov_first (f h : E → ℝ) (f' : E → E) (x0 : E) := (l : NNReal) (hl : l > (0 : ℝ)) (x y : ℕ → E) (t γ : ℕ → ℝ) (convf: ConvexOn ℝ univ f) (convh : ConvexOn ℝ univ h) (h₁ : ∀ x : E, HasGradientAt f (f' x) x) (h₂ : LipschitzWith l f')  (oriy : y 0 = x 0) (initial : x 0 = x0) (teq : ∀ n : ℕ, t n = 1 / l) (γeq : ∀ n : ℕ, γ n = 2 / (2 + n)) (update1 : ∀ k, y k = x k + (γ k * (1 - γ (k - 1)) / γ (k - 1)) • (x k - x (k - 1))) (update2 : ∀ k, prox_prop (t k • h) (y k - t k • f' (y k)) (x (k + 1)))
\end{lstlisting}
The convergence theorem is stated as:
\begin{theorem}
  Suppose that Assumption \ref{assumption} is satisfied, the fixed step size $t_k = t  \in (0, \frac{1}{L}]$,
  and the hyperparameters $\gamma_k = \frac{2}{2 + k}$, then the sequence $\{x^k\}$ generated by (\ref{alg: Nesterov}) satisfies
  \begin{align}
    \psi(x^k) - \psi^* \leq \frac{2L}{(k+1)^2} \| x^0 - x^* \|^2.
  \end{align}
\end{theorem}
We can prove the convergence rate for the Nesterov acceleration method stated by the theorem above by complicated calculation in Lean4. The formalized version of the convergence rate for the fixed stepsize is given as:

\begin{lstlisting}
theorem Nesterov_first_fix_stepsize_converge {alg: Nesterov_first f h f' x0} {xm : E} (minφ : IsMinOn (f + h) univ xm) :
    ∀ k, f (alg.x (k + 1)) + h (alg.x (k + 1)) - f alg.xm - h alg.xm ≤ 2 * alg.l / (k + 2) ^ 2 * ‖ x0 - alg.xm ‖ ^ 2
\end{lstlisting}
For the second version of Nesterov's acceleration algorithm (\ref{alg: second Nesterov}), we can formalize as:
\begin{lstlisting}
class Nesterov_second (f h : E → ℝ) (f' : E → E) (x0 : E) :=
  (l : NNReal) (hl : l > (0 : ℝ)) (x y : ℕ → E) (z : ℕ+ → E) 
  (t γ : ℕ → ℝ) (h₁ : ∀ x : E, HasGradientAt f (f' x) x) 
  (convf: ConvexOn ℝ Set.univ f) (h₂ : LipschitzWith l f')
  (convh : ConvexOn ℝ univ h) (oriy : y 0 = x 0) (oriγ : γ 1 = 1) 
  (initial : x0 = x 0) (teq : ∀ n : ℕ, t n = 1 / l)
  (γeq : ∀ n : ℕ, γ n = if n = 0 then (1 / (2 : ℝ)) else (2 : ℝ) / (1 + n))
  (update1 : ∀ k, z k = (1 - γ k) • (x (k - 1)) + γ k • (y (k - 1)))
  (update2 : ∀ k, prox_prop ((t k / γ k) • h) (y (k - 1) - (t k / γ k) • (f' (z k))) (y k)) (update3 : ∀ (k : ℕ+), x k = (1 - γ k) • (x (k - 1)) + γ k • y k)
\end{lstlisting}
For this method, we  also have the $O(\frac{1}{k^2})$ rate as:
\begin{theorem}
  Suppose that Assumption \ref{assumption} is satisfied, the fixed step size $t_k = t  \in (0, \frac{1}{L}]$,
  and the hyperparameters $\gamma_k = \frac{2}{1 + k}$, then the sequence $\{x^k\} $ generated by (\ref{alg: Nesterov}) satisfies
  \begin{align}
    \psi(x^k) - \psi^* \leq \frac{2L}{(k+1)^2} \| x^0 - x^* \|^2.
  \end{align}
\end{theorem}
The formalize version of the theorem is given as: 
\begin{lstlisting}
theorem Nesterov_second_fix_stepsize_converge {alg: Nesterov_second f h f' x0} {xm : E} (minφ : IsMinOn (f + h) univ xm): ∀ (k : ℕ), f (alg.x (k + 1)) + h (alg.x (k + 1)) - f xm - h xm ≤ 2 * alg.l / (k + 2) ^ 2 * ‖ x0 - xm ‖ ^ 2
\end{lstlisting}
\subsection{Application: Convergence Rate for Lasso Problem}
  In this subsection, we apply the formalization of the convergence rate of different algorithms for a concrete optimization problem, ``Lasso", from compressive sensing and sparse optimization. It is widely used in image processing, statistics and many other areas. Theoretical properties have been extensively studied in \cite{tibshirani2013lasso}. We demonstrate that the convergence can be easily formalized based on what we have done. The Lasso optimization problem is given as
  \begin{align}
      \min_{x \in \mathbb{R}^n} \quad  \frac{1}{2}\|Ax- b\|_2^2 + \lambda \|x\|_1,
  \end{align}
  where $A \in \mathbb{R}^{m \times n}$, $b \in \mathbb{R}^m$, and $\|x\|_1 = \sum_{i=1}^n |x_i|$ denotes the $\ell_1$-norm for the vector in $\mathbb{R}^n$. The corresponding $f(x)$ and $g(x)$ in composite optimization problem are given as $\frac{1}{2}\|Ax-b\|^2_2$  and $\lambda \|x\|_1$. From basic analysis, we can get the explicit form of $\nabla f(x) = A^T(Ax-b)$ and $\text{prox}_{g}(x)= \text{sign}(x) \max\{|x| - \lambda, 0\}$. With the explicit form of the update rule, the class of proximal gradient method for the Lasso problem can be defined as
\begin{lstlisting}
class LASSO_prox (A : Matrix (Fin m) (Fin n) ℝ) (b : (Fin m) → ℝ) 
    (μ : ℝ) (μpos : 0 < μ) (Ane0 : A ≠ 0) (x₀ : EuclideanSpace ℝ (Fin n)) :=
  (f h : EuclideanSpace ℝ (Fin n) → ℝ) 
  (x y : ℕ → (EuclideanSpace ℝ (Fin n)))
  (f' : EuclideanSpace ℝ (Fin n) → EuclideanSpace ℝ (Fin n))
  (ori : x 0 = x₀) (L : ℝ≥0) (t : ℝ) (xm : (EuclideanSpace ℝ (Fin n)))
  (feq : f = fun x => (1 / 2) * ‖ A *ᵥ x - b ‖₂ ^ 2)
  (f'eq : f' = fun x => Aᵀ *ᵥ (A *ᵥ x - b))
  (heq : h = fun y => μ • ‖ y ‖₁) (teq : t = 1 / L)
  (Leq : L = ‖(Matrix.toEuclideanLin ≪≫ₗ LinearMap.toContinuousLinearMap) (Aᵀ * A)‖₊)
  (minphi : IsMinOn (f + h) Set.univ xm) 
  (update1 : ∀ k, y k = x k - t • f' (x k))
  (update2 : ∀ k, x (k + 1) = fun i => (Real.sign (y k i) * (max (abs (y k i) - t * μ) 0)))
\end{lstlisting}
This definition contains the information of the target function, the derivative function, the Lipschitz constant and relevant update scheme. We can directly prove this update scheme for Lasso problem is exactly a special form of the proximal gradient descent method using the instance below.
\begin{lstlisting}
instance {A : Matrix (Fin m) (Fin n) ℝ} {b : (Fin m) → ℝ} {μ : ℝ} {μpos : 0 < μ} {Ane0 : A ≠ 0} {x₀ : (EuclideanSpace ℝ (Fin n))} [p : LASSO A b μ μpos Ane0 x₀] :
  proximal_gradient_method p.f p.h p.f' x₀ where
\end{lstlisting}
By the result we have in section \ref{sec: proximal gradient method}, we can easily get the convergence rate for proximal gradient method as
\begin{lstlisting}
theorem LASSO_converge {alg : LASSO A b μ μpos Ane0 x₀} :
    ∀ k, alg.f (alg.x k) + alg.h (alg.x k) - alg.f alg.xm - alg.h alg.xm ≤ alg.L / (2 * k) * ‖ x₀ - alg.xm ‖ ^ 2
\end{lstlisting}
We can easily get a similar formulation of the formalization of Nesterov's acceleration method for Lasso problem and its convergence rate with the same technique, where most of the code is the same with the class \lean{LASSO_prox} except for the particular update rules.

\section{Conclusion and Future Work}

In this paper, we primarily discuss the formalization of first-order algorithms in convex optimization. First, to conveniently demonstrate the derivative and first-order information of convex functions, we define the gradient and subgradient in Lean. Leveraging these definitions allows us to delve into the properties of convex functions and Lipschitz smooth functions. We then define the proximal operator, which is widely used in non-smooth optimization. By integrating these tools, we describe the class of first-order algorithms and prove the convergence rate for four widely used algorithms. These foundations provide a base and offer experience for proving more complex algorithms, such as ADMM \cite{wang2019global} and BCD \cite{razaviyayn2013unified}, in the near future. Future work will include defining the Fréchet sub-differentiability of general functions and the KL property, respectively. Additionally, discussing the optimality conditions of constrained optimization problems is of vital importance. We hope that, based on this work, we can progressively train the large language model to perform formalization automatically.

\section*{Acknowledgments}
Z. Wen is supported in part by the NSFC grants 12331010 and 12288101.

\section*{Author Contribution} 
C. Li, Z. Wang and Z. Wen contributed to the formalization of the algorithms, the writing and revision of the manuscript. W. He, Y. Wu and S. Xu
 contributed to the formalization of the algorithms.

\section*{Conflict Interest} The authors declare no conflict of interests.

\bibliography{sn-bibliography}


\begin{thebibliography}{26}
\ifx \bisbn   \undefined \def \bisbn  #1{ISBN #1}\fi
\ifx \binits  \undefined \def \binits#1{#1}\fi
\ifx \bauthor  \undefined \def \bauthor#1{#1}\fi
\ifx \batitle  \undefined \def \batitle#1{#1}\fi
\ifx \bjtitle  \undefined \def \bjtitle#1{#1}\fi
\ifx \bvolume  \undefined \def \bvolume#1{\textbf{#1}}\fi
\ifx \byear  \undefined \def \byear#1{#1}\fi
\ifx \bissue  \undefined \def \bissue#1{#1}\fi
\ifx \bfpage  \undefined \def \bfpage#1{#1}\fi
\ifx \blpage  \undefined \def \blpage #1{#1}\fi
\ifx \burl  \undefined \def \burl#1{\textsf{#1}}\fi
\ifx \doiurl  \undefined \def \doiurl#1{\url{https://doi.org/#1}}\fi
\ifx \betal  \undefined \def \betal{\textit{et al.}}\fi
\ifx \binstitute  \undefined \def \binstitute#1{#1}\fi
\ifx \binstitutionaled  \undefined \def \binstitutionaled#1{#1}\fi
\ifx \bctitle  \undefined \def \bctitle#1{#1}\fi
\ifx \beditor  \undefined \def \beditor#1{#1}\fi
\ifx \bpublisher  \undefined \def \bpublisher#1{#1}\fi
\ifx \bbtitle  \undefined \def \bbtitle#1{#1}\fi
\ifx \bedition  \undefined \def \bedition#1{#1}\fi
\ifx \bseriesno  \undefined \def \bseriesno#1{#1}\fi
\ifx \blocation  \undefined \def \blocation#1{#1}\fi
\ifx \bsertitle  \undefined \def \bsertitle#1{#1}\fi
\ifx \bsnm \undefined \def \bsnm#1{#1}\fi
\ifx \bsuffix \undefined \def \bsuffix#1{#1}\fi
\ifx \bparticle \undefined \def \bparticle#1{#1}\fi
\ifx \barticle \undefined \def \barticle#1{#1}\fi
\bibcommenthead
\ifx \bconfdate \undefined \def \bconfdate #1{#1}\fi
\ifx \botherref \undefined \def \botherref #1{#1}\fi
\ifx \url \undefined \def \url#1{\textsf{#1}}\fi
\ifx \bchapter \undefined \def \bchapter#1{#1}\fi
\ifx \bbook \undefined \def \bbook#1{#1}\fi
\ifx \bcomment \undefined \def \bcomment#1{#1}\fi
\ifx \oauthor \undefined \def \oauthor#1{#1}\fi
\ifx \citeauthoryear \undefined \def \citeauthoryear#1{#1}\fi
\ifx \endbibitem  \undefined \def \endbibitem {}\fi
\ifx \bconflocation  \undefined \def \bconflocation#1{#1}\fi
\ifx \arxivurl  \undefined \def \arxivurl#1{\textsf{#1}}\fi
\csname PreBibitemsHook\endcsname

\bibitem[\protect\citeauthoryear{Bottou et~al.}{2018}]{bottou2018optimization}
\begin{barticle}
\bauthor{\bsnm{Bottou}, \binits{L.}},
\bauthor{\bsnm{Curtis}, \binits{F.E.}},
\bauthor{\bsnm{Nocedal}, \binits{J.}}:
\batitle{Optimization methods for large-scale machine learning}.
\bjtitle{SIAM review}
\bvolume{60}(\bissue{2}),
\bfpage{223}--\blpage{311}
(\byear{2018})
\end{barticle}
\endbibitem

\bibitem[\protect\citeauthoryear{Boldo et~al.}{2016}]{boldo2016formalization}
\begin{barticle}
\bauthor{\bsnm{Boldo}, \binits{S.}},
\bauthor{\bsnm{Lelay}, \binits{C.}},
\bauthor{\bsnm{Melquiond}, \binits{G.}}:
\batitle{Formalization of real analysis: A survey of proof assistants and libraries}.
\bjtitle{Mathematical Structures in Computer Science}
\bvolume{26}(\bissue{7}),
\bfpage{1196}--\blpage{1233}
(\byear{2016})
\end{barticle}
\endbibitem

\bibitem[\protect\citeauthoryear{Nipkow et~al.}{2002}]{Nipkow2002APA}
\begin{botherref}
\oauthor{\bsnm{Nipkow}, \binits{T.}},
\oauthor{\bsnm{Paulson}, \binits{L.C.}},
\oauthor{\bsnm{Wenzel}, \binits{M.}}:
A proof assistant for higher-order logic.
Lecture Notes in Computer Science
(2002)
\end{botherref}
\endbibitem

\bibitem[\protect\citeauthoryear{De~Moura et~al.}{2015}]{de2015lean}
\begin{bchapter}
\bauthor{\bsnm{De~Moura}, \binits{L.}},
\bauthor{\bsnm{Kong}, \binits{S.}},
\bauthor{\bsnm{Avigad}, \binits{J.}},
\bauthor{\bsnm{Van~Doorn}, \binits{F.}},
\bauthor{\bsnm{Raumer}, \binits{J.}}:
\bctitle{The lean theorem prover (system description)}.
In: \bbtitle{Automated Deduction-CADE-25: 25th International Conference on Automated Deduction, Berlin, Germany, August 1-7, 2015, Proceedings 25},
pp. \bfpage{378}--\blpage{388}
(\byear{2015}).
\bcomment{Springer}
\end{bchapter}
\endbibitem

\bibitem[\protect\citeauthoryear{Kudryashov}{2022}]{kudryashov2022formalizing}
\begin{bchapter}
\bauthor{\bsnm{Kudryashov}, \binits{Y.}}:
\bctitle{Formalizing the divergence theorem and the cauchy integral formula in lean}.
In: \bbtitle{13th International Conference on Interactive Theorem Proving (ITP 2022)}
(\byear{2022}).
\bcomment{Schloss-Dagstuhl-Leibniz Zentrum f{\"u}r Informatik}
\end{bchapter}
\endbibitem

\bibitem[\protect\citeauthoryear{Gou{\"e}zel}{2022}]{gouezel2022formalization}
\begin{bchapter}
\bauthor{\bsnm{Gou{\"e}zel}, \binits{S.}}:
\bctitle{A formalization of the change of variables formula for integrals in mathlib}.
In: \bbtitle{International Conference on Intelligent Computer Mathematics},
pp. \bfpage{3}--\blpage{18}
(\byear{2022}).
\bcomment{Springer}
\end{bchapter}
\endbibitem

\bibitem[\protect\citeauthoryear{Tassarotti et~al.}{2021}]{tassarotti2021formal}
\begin{bchapter}
\bauthor{\bsnm{Tassarotti}, \binits{J.}},
\bauthor{\bsnm{Vajjha}, \binits{K.}},
\bauthor{\bsnm{Banerjee}, \binits{A.}},
\bauthor{\bsnm{Tristan}, \binits{J.-B.}}:
\bctitle{A formal proof of pac learnability for decision stumps}.
In: \bbtitle{Proceedings of the 10th ACM SIGPLAN International Conference on Certified Programs and Proofs},
pp. \bfpage{5}--\blpage{17}
(\byear{2021})
\end{bchapter}
\endbibitem

\bibitem[\protect\citeauthoryear{Grechuk}{2011}]{Lower_Semicontinuous-AFP}
\begin{botherref}
\oauthor{\bsnm{Grechuk}, \binits{B.}}:
Lower semicontinuous functions.
Archive of Formal Proofs
(2011).
\url{https://isa-afp.org/entries/Lower_Semicontinuous.html}, Formal proof development
\end{botherref}
\endbibitem

\bibitem[\protect\citeauthoryear{Allamigeon and Katz}{2019}]{allamigeon2019formalization}
\begin{barticle}
\bauthor{\bsnm{Allamigeon}, \binits{X.}},
\bauthor{\bsnm{Katz}, \binits{R.D.}}:
\batitle{A formalization of convex polyhedra based on the simplex method}.
\bjtitle{J. Automat. Reason.}
\bvolume{63}(\bissue{2}),
\bfpage{323}--\blpage{345}
(\byear{2019})
\doiurl{10.1007/s10817-018-9477-1}
\end{barticle}
\endbibitem

\bibitem[\protect\citeauthoryear{Bentkamp et~al.}{2023}]{verifiedoptimization}
\begin{bchapter}
\bauthor{\bsnm{Bentkamp}, \binits{A.}},
\bauthor{\bsnm{Mir}, \binits{R.F.}},
\bauthor{\bsnm{Avigad}, \binits{J.}}:
\bctitle{Verified reductions for optimization}.
In: \beditor{\bsnm{Sankaranarayanan}, \binits{S.}},
\beditor{\bsnm{Sharygina}, \binits{N.}} (eds.)
\bbtitle{Tools and Algorithms for the Construction and Analysis of Systems},
pp. \bfpage{74}--\blpage{92}.
\bpublisher{Springer},
\blocation{Cham}
(\byear{2023})
\end{bchapter}
\endbibitem

\bibitem[\protect\citeauthoryear{mathlib Community}{2020}]{mathlibcommunity}
\begin{bchapter}
\bauthor{\bsnm{Community}, \binits{T.}}:
\bctitle{The lean mathematical library}.
In: \bbtitle{Proceedings of the 9th ACM SIGPLAN International Conference on Certified Programs and Proofs}.
\bsertitle{CPP 2020},
pp. \bfpage{367}--\blpage{381}.
\bpublisher{Association for Computing Machinery},
\blocation{New York, NY, USA}
(\byear{2020}).
\doiurl{10.1145/3372885.3373824} .
\burl{https://doi.org/10.1145/3372885.3373824}
\end{bchapter}
\endbibitem

\bibitem[\protect\citeauthoryear{Nesterov}{1983}]{nesterov1983method}
\begin{barticle}
\bauthor{\bsnm{Nesterov}, \binits{Y.}}:
\batitle{A method of solving a convex programming problem with convergence rate o (1/k** 2)}.
\bjtitle{Doklady Akademii Nauk SSSR}
\bvolume{269}(\bissue{3}),
\bfpage{543}
(\byear{1983})
\end{barticle}
\endbibitem

\bibitem[\protect\citeauthoryear{Liu et~al.}{2020}]{optimization_method}
\begin{bbook}
\bauthor{\bsnm{Liu}, \binits{H.}},
\bauthor{\bsnm{Hu}, \binits{J.}},
\bauthor{\bsnm{Li}, \binits{Y.}},
\bauthor{\bsnm{Wen}, \binits{Z.}}:
\bbtitle{Optimization: Modeling, Algorithm and Theory}.
\bpublisher{higher education press},
\blocation{Beijing, China}
(\byear{2020})
\end{bbook}
\endbibitem

\bibitem[\protect\citeauthoryear{Nesterov}{2018}]{nesterov2018lectures}
\begin{bbook}
\bauthor{\bsnm{Nesterov}, \binits{Y.}}:
\bbtitle{Lectures on Convex Optimization},
\bedition{2nd} edn.
\bpublisher{Springer},
\blocation{Gewerbestrasse 11, 6330 Cham, Switzerland}
(\byear{2018})
\end{bbook}
\endbibitem

\bibitem[\protect\citeauthoryear{Nocedal and Wright}{1999}]{nocedal1999numerical}
\begin{bbook}
\bauthor{\bsnm{Nocedal}, \binits{J.}},
\bauthor{\bsnm{Wright}, \binits{S.J.}}:
\bbtitle{Numerical Optimization},
\bedition{2nd} edn.
\bpublisher{Springer},
\blocation{New York, NY 10013, USA}
(\byear{1999})
\end{bbook}
\endbibitem

\bibitem[\protect\citeauthoryear{Beck}{2017}]{beck2017first}
\begin{bbook}
\bauthor{\bsnm{Beck}, \binits{A.}}:
\bbtitle{First-Order Methods in Optimization}.
\bpublisher{SIAM-Society for Industrial and Applied Mathematics},
\blocation{Philadelphia, PA, USA}
(\byear{2017})
\end{bbook}
\endbibitem

\bibitem[\protect\citeauthoryear{Sun and Yuan}{2006}]{sun2006optimization}
\begin{bbook}
\bauthor{\bsnm{Sun}, \binits{W.}},
\bauthor{\bsnm{Yuan}, \binits{Y.-X.}}:
\bbtitle{Optimization Theory and Methods: Nonlinear Programming}
vol. \bseriesno{1}.
\bpublisher{Springer},
\blocation{New York, NY 10013, USA}
(\byear{2006})
\end{bbook}
\endbibitem

\bibitem[\protect\citeauthoryear{Boyd and Vandenberghe}{2004}]{boyd2004convex}
\begin{bbook}
\bauthor{\bsnm{Boyd}, \binits{S.P.}},
\bauthor{\bsnm{Vandenberghe}, \binits{L.}}:
\bbtitle{Convex Optimization}.
\bpublisher{Cambridge university press},
\blocation{Cambridge, CB2 8RU, UK}
(\byear{2004})
\end{bbook}
\endbibitem

\bibitem[\protect\citeauthoryear{Baanen}{2022}]{baanen2022use}
\begin{bchapter}
\bauthor{\bsnm{Baanen}, \binits{A.}}:
\bctitle{{Use and Abuse of Instance Parameters in the Lean Mathematical Library}}.
In: \beditor{\bsnm{Andronick}, \binits{J.}},
\beditor{\bsnm{Moura}, \binits{L.}} (eds.)
\bbtitle{13th International Conference on Interactive Theorem Proving (ITP 2022)}.
\bsertitle{Leibniz International Proceedings in Informatics (LIPIcs)},
vol. \bseriesno{237},
pp. \bfpage{4}--\blpage{1420}.
\bpublisher{Schloss Dagstuhl -- Leibniz-Zentrum f{\"u}r Informatik},
\blocation{Dagstuhl, Germany}
(\byear{2022}).
\doiurl{10.4230/LIPIcs.ITP.2022.4} .
\burl{https://drops.dagstuhl.de/entities/document/10.4230/LIPIcs.ITP.2022.4}
\end{bchapter}
\endbibitem

\bibitem[\protect\citeauthoryear{Boyd et~al.}{2003}]{boyd2003subgradient}
\begin{botherref}
\oauthor{\bsnm{Boyd}, \binits{S.}},
\oauthor{\bsnm{Xiao}, \binits{L.}},
\oauthor{\bsnm{Mutapcic}, \binits{A.}}:
Subgradient methods.
lecture notes of EE392o, Stanford University, Autumn Quarter
\textbf{2004}(01)
(2003)
\end{botherref}
\endbibitem

\bibitem[\protect\citeauthoryear{Parikh et~al.}{2014}]{parikh2014proximal}
\begin{barticle}
\bauthor{\bsnm{Parikh}, \binits{N.}},
\bauthor{\bsnm{Boyd}, \binits{S.}}, \betal:
\batitle{Proximal algorithms}.
\bjtitle{Foundations and trends{\textregistered} in Optimization}
\bvolume{1}(\bissue{3}),
\bfpage{127}--\blpage{239}
(\byear{2014})
\end{barticle}
\endbibitem

\bibitem[\protect\citeauthoryear{Beck and Teboulle}{2009}]{beck2009fast}
\begin{barticle}
\bauthor{\bsnm{Beck}, \binits{A.}},
\bauthor{\bsnm{Teboulle}, \binits{M.}}:
\batitle{A fast iterative shrinkage-thresholding algorithm for linear inverse problems}.
\bjtitle{SIAM journal on imaging sciences}
\bvolume{2}(\bissue{1}),
\bfpage{183}--\blpage{202}
(\byear{2009})
\end{barticle}
\endbibitem

\bibitem[\protect\citeauthoryear{Bolte et~al.}{2014}]{bolte2014proximal}
\begin{barticle}
\bauthor{\bsnm{Bolte}, \binits{J.}},
\bauthor{\bsnm{Sabach}, \binits{S.}},
\bauthor{\bsnm{Teboulle}, \binits{M.}}:
\batitle{Proximal alternating linearized minimization for nonconvex and nonsmooth problems}.
\bjtitle{Mathematical Programming}
\bvolume{146}(\bissue{1-2}),
\bfpage{459}--\blpage{494}
(\byear{2014})
\end{barticle}
\endbibitem

\bibitem[\protect\citeauthoryear{Tibshirani}{2013}]{tibshirani2013lasso}
\begin{barticle}
\bauthor{\bsnm{Tibshirani}, \binits{R.J.}}:
\batitle{{The lasso problem and uniqueness}}.
\bjtitle{Electronic Journal of Statistics}
\bvolume{7}(\bissue{none}),
\bfpage{1456}--\blpage{1490}
(\byear{2013})
\doiurl{10.1214/13-EJS815}
\end{barticle}
\endbibitem

\bibitem[\protect\citeauthoryear{Wang et~al.}{2019}]{wang2019global}
\begin{barticle}
\bauthor{\bsnm{Wang}, \binits{Y.}},
\bauthor{\bsnm{Yin}, \binits{W.}},
\bauthor{\bsnm{Zeng}, \binits{J.}}:
\batitle{Global convergence of admm in nonconvex nonsmooth optimization}.
\bjtitle{Journal of Scientific Computing}
\bvolume{78},
\bfpage{29}--\blpage{63}
(\byear{2019})
\end{barticle}
\endbibitem

\bibitem[\protect\citeauthoryear{Razaviyayn et~al.}{2013}]{razaviyayn2013unified}
\begin{barticle}
\bauthor{\bsnm{Razaviyayn}, \binits{M.}},
\bauthor{\bsnm{Hong}, \binits{M.}},
\bauthor{\bsnm{Luo}, \binits{Z.-Q.}}:
\batitle{A unified convergence analysis of block successive minimization methods for nonsmooth optimization}.
\bjtitle{SIAM J. Optim.}
\bvolume{23}(\bissue{2}),
\bfpage{1126}--\blpage{1153}
(\byear{2013})
\doiurl{10.1137/120891009}
\end{barticle}
\endbibitem

\end{thebibliography}

\end{document}